\documentclass[12pt]{amsart}
\usepackage{amssymb}
\input amssym.def
\usepackage{amsmath, amsfonts,hyperref,xcolor}
\usepackage{relsize}
\usepackage{amscd}
\usepackage[mathscr]{eucal}
\usepackage{orcidlink}
\hypersetup{colorlinks=true,citecolor={purple},linkcolor={blue},urlcolor={violet}}

\makeatletter
\@namedef{subjclassname@2020}{%
  \textup{2020} Mathematics Subject Classification}
\makeatother

\newcommand{\N}{{\mathbb N}}
\newcommand{\Z}{{\mathbb Z}}

\newtheorem{thm}{Theorem}
\newtheorem{lem}{Lemma}
\newtheorem{cor}{Corollary}

\newtheorem{rmk}{Remark}
\newtheorem{defn}{Definition}
\newcommand{\thmref}[1]{Theorem~\ref{#1}}

\newcommand{\lemref}[1]{Lemma~\ref{#1}}
\newcommand{\corref}[1]{Corollary~\ref{#1}}

\begin{document}

\title[Additive convolution sum of the arithmetic function]{On additive convolution sum of arithmetic functions and related questions}

\author[B. Misra, B. Saha and A. Sharma]{Bikram Misra \orcidlink{0009-0000-5863-2789}, Biswajyoti Saha \orcidlink{0009-0009-2904-4860} and Anubhav Sharma
\orcidlink{0000-0002-1522-4432}}

\address{Department of Mathematics, Indian Institute of Technology Delhi, New Delhi 110016, India}
\email{maz218517@iitd.ac.in, biswajyoti@maths.iitd.ac.in}

\address{Theoretical Statistics and Mathematics Unit, Indian Statistical Institute, Delhi Center, New Delhi 110016, India}
\email{anubhav6595@gmail.com}



\maketitle

\begin{abstract}
Ingham studied two types of convolution sums of the divisor function, 
the shifted convolution sum $\sum_{n \le N} d(n) d(n+h)$ and the additive convolution sum
$\sum_{n < N} d(n) d(N-n)$ for integers $N, h$ and derived their asymptotic formulas as $N \to \infty$.
There have been numerous works extending Ingham's result on the shifted convolution sum,
but only little has been done towards
the additive convolution sum. In this article, we extend the classical result of Ingham 
to derive an asymptotic formula with an error term of the sub-sum $\sum_{n < M} d(n) d(N-n)$ for certain integers $M \le N$. This involves careful choice of an applicable range of $M$. 
We also study the convolution sum $\sum_{n < M} f(n) g(N-n)$ for certain arithmetic functions $f$ and $g$
with absolutely convergent Ramanujan expansions, which in turn leads us to a well-established prediction of Ramanujan.
\end{abstract}

\section{Introduction}
Convolution sums are ubiquitous in number theory and they play very crucial role
in numerous arithmetic problems. Two such instances come from
two of the long-standing open problems in number theory.
In the context of the twin prime conjecture, one
considers the shifted convolution sum
$$
\sum_{n \le N} \Lambda(n) \Lambda(n+h),
$$
where $\Lambda(n)$ denotes the von Mangoldt function, $N, h$ are positive integers, with $h$ even. Note that in the above sum,
the summand is zero unless both $n$ and $n+h$ are prime powers and the contribution of the prime
powers is small compared to the contribution of the primes themselves. Hence an asymptotic estimate
of this sum would enable us to deduce the infinitude of the prime pairs of the form $(n,n+h)$,
giving us the  twin prime conjecture by putting $h=2$.

The next instance is in the context of the Goldbach conjecture. Here one considers the additive convolution sum
$$
\sum_{n < N} \Lambda(n) \Lambda(N-n),
$$
for an even positive integer $N$. Note that this sum counts the number of ways $N$ can be written as a sum of two
prime powers. So in this case, an asymptotic estimate of the above sum leads to  the Goldbach conjecture.
Moreover, here it is also relevant to consider the sub-sum 
$$
\sum_{n \le M} \Lambda(n) \Lambda(N-n),
$$
where $M \le N/2$. To elaborate, let $N \ge 4$ be a positive even integer which can be written
as a sum of two primes $p_1,p_2$. We refer to $p_1, p_2$ as prime summands of $N$. Clearly, $\min\{p_1,p_2\} \le N/2$.
So assuming the Goldbach conjecture, we are assured to get a prime summand of $N$ which
is $\le N/2$. Note that it may be possible to write $N$ as a sum of two primes in many ways,
for example, $14=3+11=7+7$. Hence for determining the least possible prime summand of $N$,
one can consider the sum $\sum_{n \le M} \Lambda(n) \Lambda(N-n)$ and study
its behaviour with respect to both $M$ and $N$, where $M$ is a parameter $\le N/2$. The study of
such convolution sub-sums can also be motivated by the following theorem of Chen \cite{JRC} that every sufficiently
large even number can be written as a sum of a prime and a semi-prime. By a semi-prime we mean a positive integer
which is product of at most two primes. Chen's result
was further strengthened by Cai \cite{YCC} to show that every sufficiently large even integer $n$ is a sum of a prime
$\le n^{0.95}$ and a semi-prime.

Also, if we are to study the sums of the form $\sum_{n <N} f(n) g(N-n) h(n)$ for some suitable arithmetic functions
$f,g$ and $h$, then we can estimate such a  sum using the partial summation formula if we have information about
the sub-sum $\sum_{n <M} f(n) g(N-n)$ for $M \le N$.
Hence in this article, we take up the study of the additive convolution sub-sum $\sum_{n <M} f(n) g(N-n)$
of arithmetic functions $f,g$ for integers $M \le N$. We begin with the divisor function and the classical result of Ingham \cite{AEI}.
In the seminal paper \cite{AEI}, Ingham studied the following convolution sums of the divisor function
$d(n)=\sum_{d \mid n \atop d>0} 1$. He proved the following asymptotic formulas:
\begin{equation}\label{shifted-divisor}
\sum_{n \le N} d(n) d(n+h) = \frac{6}{\pi^2} \sigma_{-1}(h) N \log^2 N+ O(N \log N),
\end{equation}
for $N \to \infty$ and 
\begin{equation}\label{add-divisor-1}
\sum_{n < N} d(n) d(N-n) = \frac{6}{\pi^2} \sigma_{1}(N) \log^2 N + O(\sigma_1(N) \log N \log\log N),
\end{equation}
for $N \to \infty$, where $\displaystyle{\sigma_s(n):=\sum_{d \mid n\atop d>0} d^s}$ for a complex number $s$.

This fundamental work \cite{AEI} has been a central force of developments of many key concepts
in the study of the mean values of some important $L$-functions (for example \cite{DI, YM}) and
further led to many important studies of critical values, shifted convolution sums and subconvexity questions.
However, as Motohashi \cite{YM} writes, the research on the additive convolution sums has not been
as intensive as on the shifted convolution sums, perhaps because of no apparent relation
with the mean values of $L$-functions. A recent work to be noted here is \cite{AM}.
Estermann \cite{TE} provided the secondary main terms in \eqref{add-divisor-1}, using methods different
from \cite{AEI}.

In this article, we particularly consider the last formula \eqref{add-divisor-1} and derive an extension. Note that,
except perhaps one summand of $\sum_{n < N} d(n) d(N-n)$
in the middle (which has insignificant contribution), all the summands appear twice and
due to the palindromic nature of the sum, the asymptotic formula \eqref{add-divisor-1} is equivalent to
\begin{equation}\label{add-divisor-2}
\sum_{n \le N/2} d(n) d(N-n) = \frac{3}{\pi^2} \sigma_{1}(N) \log^2 N + O(\sigma_1(N) \log N \log\log N),
\end{equation}
for $N \to \infty$. Taking motivation from our earlier discussions, we find it interesting to consider the sub-sum
$\sum_{n \le M} d(n) d(N-n)$ for an integer $M \le N$. It might be natural to expect the contribution of the
sub-sum to be $M/N$ times that of the complete sum. In addition, we also get an
appropriate change in the logarithmic factor in the main term,
but more importantly, we need to choose the range of $M$ carefully.

Following Ingham's method \cite{AEI}, we first derive an
asymptotic formula of $\sum_{n \le M} d(n) d(N-n)$ with an error term which depends on $N$ and $M$, for $N$ sufficiently large and $M\le N/2$ such that ${N}/{(\sigma_{-1}(N)\log N)^2}=o(M)$.
More precisely, we prove:

\begin{thm}\label{thm-add-divisor}
For $N$ sufficiently large and $M\le N/2$ such that ${N}/{(\sigma_{-1}(N)\log N)^2}=o(M)$, we have
\begin{equation}
\sum_{n\le M}d(n)d(N-n) \sim \frac{6}{\pi^2}M\sigma_{-1}(N)\log ^2X,
\end{equation}
where $X=\sqrt{M(N-M)}$.
Moreover, for ${N}/{(\sigma_{-1}(N)\log\log N)^2}\le M\le {N}/{2}$, we have
\begin{equation}\label{ex-add-divisor-1}
\sum_{n\le M} d(n)d(N-n) =\frac{6}{\pi^2}M\sigma_{-1}(N)\log ^2X +O(M\sigma_{-1}(N)\log N \log\log N),
\end{equation}
and for $M\le {N}/{(\sigma_{-1}(N)\log\log N)^2}$, such that ${N}/{(\sigma_{-1}(N)\log N)^2}=o(M)$ we have
\begin{equation}\label{ex-add-divisor-2}
\sum_{n\le M} d(n)d(N-n) =\frac{6}{\pi^2}M\sigma_{-1}(N)\log ^2X +O(X\log X).
\end{equation}
\end{thm}

\begin{rmk}\rm
If we put $M=N/2$ in Theorem \ref{thm-add-divisor}, we get back Ingham's result \eqref{add-divisor-2}.
\end{rmk}

Note that, for any real number $M$ such that $1 \le M \le N$, we have
$$
\sum_{n< M} d(n)d(N-n)
=\sum_{n< N} d(n)d(N-n)-\sum_{n \le N-M}d(n)d(N-n).
$$
Now if $N/2 \le M \le N$, we have $N-M \le N/2$.  So as a corollary to \thmref{thm-add-divisor}, we get the following:

\begin{cor}\label{cor-add-divisor}
For sufficiently large $M$ such that $N/2 \le M \le N$ and ${N}/{(\sigma_{-1}(N)\log N)^2}=o(N-M)$, we have
\begin{equation}\label{ex-add-divisor-2}
\sum_{n< M}d(n)d(N-n) \sim\frac{6}{\pi^2}\sigma_{-1}(N)\left(N\log^2N-(N-M)\log^2 X\right),
\end{equation}
where $X=\sqrt{M(N-M)}$. Moreover, for $N/2\le M\le N-{N}/{(\sigma_{-1}(N)\log\log N)^2}$, we have
\begin{equation}
\sum_{n< M}d(n)d(N-n) =\frac{6}{\pi^2}\sigma_{-1}(N)\left(N\log^2N-(N-M)\log^2 X\right)+O(\sigma_{1}(N)\log N\log\log N).
\end{equation}
and for $M\ge N-{N}/{(\sigma_{-1}(N)\log\log N)^2}$ such that ${N}/{\sigma_{-1}(N)^2\log^2 N}=o(N-M)$ we have
\begin{equation}
\sum_{n< M}d(n)d(N-n) =\frac{6}{\pi^2}\sigma_{-1}(N)\left(N\log^2N-(N-M)\log^2 X\right)+O(X\log X).
\end{equation}
\end{cor}

The above statements give key insights regarding the sub-sum $\sum_{n \le M}d(n)d(N-n)$ in comparison with the 
full-sum $\sum_{n < N} d(n) d(N-n)$. For example, if $M \sim c N$ for some $0 <c <1$, then asymptotically the sub-sum is
$c$ times the full-sum. Although the above theorem includes some choices of $M$ that are $o(N)$,
the proof however seems to show that this method cannot treat the small values of $M$ of the form $N^\alpha$ for some $\alpha <1$.

\subsection{$M=N^\alpha$ for some $\alpha <1$}
This is a delicate problem and one can proceed as follows. We write
\begin{align*}
\sum_{n \le M}d(n)d(N-n)
&= \sum_{n \le M}d(n) \left( 2 \sum_{t|N-n \atop t \le \sqrt{N-n}} 1 +O(1)\right)\\
&= 2 \sum_{n \le M}d(n) \left( \sum_{t|N-n \atop t \le \sqrt{N}} 1 -\sum_{t|N-n \atop \sqrt{N-n} < t \le \sqrt{N}} 1\right) +O(M \log M).
\end{align*}
It is not difficult to see that 
\begin{equation}\label{small-part}
    \sum_{n \le M}d(n)\sum_{t|N-n \atop \sqrt{N-n} < t \le \sqrt{N}}1 \ll M \log M.
\end{equation}
Hence we need to study
\begin{equation}\label{eq-MT}
2 \sum_{n \le M}d(n) \sum_{t|N-n \atop t \le \sqrt{N}} 1
=2 \sum_{t \le \sqrt{N}} \sum_{n \le M \atop n \equiv N \bmod t} d(n).
\end{equation}
This therefore requires us to estimate the divisor function on residue classes.
Completing an estimate of Fouvry and Iwaniec, in \cite{PV} Pongsriiam and Vaughan  proved that
\begin{equation}\label{eq-PV}
\sum_{n \le M \atop n \equiv N \bmod t} d(n)
= \frac{M}{t} \sum_{r \mid t} \frac{c_r(N)}{r} \left( \log\frac{M}{r^2} + 2\gamma-1 \right) + O((M^{1/3}+t^{1/2})M^{\epsilon}),
\end{equation}
uniformly for $t < M^{2/3 - \epsilon}$ and arbitrarily small $\epsilon$, where $\gamma$ is the Euler's constant and $c_r(N)$ is the Ramanujan sum $\sum_{a\in (\Z/r\Z)^*}e^{2\pi ian/r}$. Now to apply \eqref{eq-PV} in \eqref{eq-MT}, we need $\sqrt N < M^{2/3 - \epsilon}$, equivalently, $N^{3/4+\epsilon}<M$. Hence, for $3/4 < \alpha <1$, we can therefore treat the additive convolution sub-sum $\sum_{n \le M}d(n)d(N-n)$ and prove the following theorem.

\begin{thm}\label{thm-add-divisor-2}
For $N$ sufficiently large and $M=N^\alpha$ such that $3/4 < \alpha <1$, we have
\begin{equation*}
\sum_{n\le M}d(n)d(N-n) = M \alpha (\log N)^2 \mathfrak S_0(N)
-2M (1+\alpha)\log N \mathfrak S_1(N) + 4M \mathfrak S_2(N)
+ O(M \log M),
\end{equation*}
where $\mathfrak S_i(N)$ is the series $\sum_{d \le \sqrt N} c_d(N) (\log d)^i/d^2$ for $i=0,1,2$, and the implied constant depends on $N$.
\end{thm}


\section{Additive convolution sums of general arithmetic functions}

Next we consider such additive convolution sums for general arithmetic functions. To make a meaningful statement, we 
need some conditions on these functions. Here we consider
arithmetic functions with absolutely convergent Ramanujan expansions. In \cite{SR}, Ramanujan studied the
Ramanujan sums introduced above and initiated the theory of Ramanujan expansions.
For positive integers $r,n$, Ramanujan sum is defined by
$$
c_r(n):=\sum_{a\in (\Z/r\Z)^*}\zeta_r^{an},
$$
where $\zeta_r$ denotes a primitive $r$-th root of unity. Since then, these sums are called Ramanujan sums. It is possible to express $c_r(n)$ in
terms of the M\"obius function $\mu$ as follows:
\begin{equation}\label{crn-exp1}
c_r(n)=\sum_{d \mid \gcd(n,r)} \mu\left(\frac{r}{d}\right) d.
\end{equation}
Ramanujan studied these sums to derive point-wise convergent series expansion of the form
$\sum_{r\ge 1} a_rc_r(n)$ for various arithmetic functions.
These expansions are now known as Ramanujan expansions. More precisely:

\begin{defn}
Let $f$ be an arithmetic function. Then $f$ is said to
admit a Ramanujan expansion if for each $n$,
$f(n)$ can be written as a convergent series
of the form
$$
f(n)= \sum_{r \ge 1} a_f(r) c_r(n)
$$
for appropriate complex numbers $a_f(r)$.
The number $a_f(r)$ is said to be the $r$-th Ramanujan coefficient
of $f$ with respect to this expansion.
\end{defn}

Ramanujan \cite{SR} showed that, for a real variable $s > 0$,
\begin{equation}\label{sigma_s}
 \frac{\sigma_s(n)}{n^s}=\zeta(s+1)\sum_{r\ge 1} \frac{c_r(n)}{r^{s+1}}.
\end{equation}
Here $\zeta(s)$ denotes the Riemann zeta function. Also, for the divisor function $d(n)$, one has
$$
d(n)=-\sum_{r \ge 1} \frac{\log r}{r} c_r(n).
$$
Ramanujan \cite{SR} himself observed that such an expansion is not necessarily unique. For example, we have
$$
\sum_{r \ge 1} \frac{c_r(n)}{r} = 0,
$$
and this equation can be viewed as a Ramanujan expansion of the zero function.
For an exposition on Ramanujan sums and Ramanujan expansions, we refer to \cite{RM}.


A systematic study of the shifted convolution sum $\sum_{n\leq N} f(n) g(n+h)$
for certain arithmetic functions $f$ and $g$ with absolutely convergent Ramanujan expansions
has been carried out in \cite{GMP,MS,CMS,BS}. However, the study of the additive convolution sum
$\sum_{n \le M} f(n) g(N-n)$ for a real number $M < N$ and such arithmetic functions $f$ and $g$
has not been carried out so far with great detail, except for an indicative statement in \cite[Theorem 10]{GMP}.
Now we derive an asymptotic formula with an error term
for the additive convolution sum $\sum_{n < M} f(n) g(N-n)$ for $M=o( N)$ and certain arithmetic functions $f$ and $g$ with absolutely convergent Ramanujan expansions with Ramanujan coefficients having appropriate growth conditions. More precisely, we prove the following:

\begin{thm}\label{thm-add-conv}
Let $f$ and $g$ be two arithmetic functions with absolutely convergent Ramanujan expansions:
$$
f(n)=\sum_{r \ge 1}a_{f}(r)c_r(n), \quad g(n)=\sum_{s \ge 1} a_{g}(s)c_s(n).
$$
Suppose that $|a_{f}(r)|, |a_{g}(r)| \ll  r^{-(1+\delta)}$ for some $\delta>0$. For $M < N$
sufficiently large and $M=o(N)$, we have
$$\sum_{n< M}f(n)g(N-n)=\begin{cases}
M\sum_{r \ge 1}a_{f}(r)a_{g}(r)c_r(N)+O(M^{1-\delta}(\log M)^{4-2\delta}) & \text{if $\delta<1$},\\
M\sum_{r \ge 1}a_{f}(r)a_{g}(r)c_r(N)+ O(\log^3M)& \text{if $\delta=1$}, \\
M\sum_{r \ge 1}a_{f}(r)a_{g}(r)c_r(N)+ O(1)& \text{if $\delta>1$},
\end{cases}
$$
where the implied constants depend on $N$.
\end{thm}

In the following section we first present a proof of Theorem \ref{thm-add-divisor}. Then in Section 4 and 5, we give a proof of \thmref{thm-add-divisor-2} and \thmref{thm-add-conv}, respectively.

\section{Proof of Theorem \ref{thm-add-divisor}}

For the proof of Theorem \ref{thm-add-divisor}, we follow Ingham's idea, but we need to make
important changes as we go ahead. The proof also highlights why it has been important to restrict
our statement for $M \le N/2$.
We write the sum 
$$
\sum_{n\leq M} d(n)d(N-n)=\sum_{n\leq M} \left(\sum_{m\mid n} 1\right)\left(\sum_{l\mid N-n} 1\right).
$$
So this is same as counting the number of $4$-tuples $(l,m,r,s)\in \N^4$ satisfying $lr+ms=N$
such that $n\leq M$, where $n=ms$ and $N-n=lr$.
Hence our goal is to find the cardinality of the set 
$$S=\{(l,r,m,s)\in \N^4 \mid lr+ms=N; \ ms\leq M\}.$$
Now, $lm \cdot rs=ms \cdot lr=n(N-n)$.
Since $n \le M$, we get $N - n \ge N - M$. Note that this does not offer any upper bound for the product $n(N-n)$.
Using AM-GM inequality, Ingham \cite{AEI} obtained the upper bound $N^2/4$. Since we have the restriction
$n \le M$, this upper bound is not optimal in terms of $M$.
So we consider the function $f(x):=x(N-x)$ on the interval $[0,N]$. This function is increasing on the sub-interval $[0,N/2]$
(and hence on $[0,M]$ as $M \le N/2$). We therefore get the bound
$$
lm \cdot rs=ms \cdot lr=n(N-n) \le M (N-M).
$$
For $n=M$ (when $M$ is an integer), this upper bound is achieved, and hence this bound is optimal.
Denote $X= \sqrt{ M (N-M)}$. Since $lm \cdot rs \le X^2$, 
at least one of the inequalities $lm\leq X$ and $rs\leq X$ must be satisfied.
So, we have $|S|=|P|+|Q|-|R|$, where
$$
P=\{(l,r,m,s)\in \N^4 \mid lr+ms=N; \ ms\leq M; \ lm\leq X\},
$$
$$
Q=\{(l,r,m,s)\in \N^4 \mid lr+ms=N; \ ms\leq M; \ rs\leq X\},
$$
and
$$
R=\{(l,r,m,s)\in \N^4 \mid lr+ms=N; \ ms\leq M; \ lm, rs\leq X\}.
$$

Since by symmetry, we have $|Q|=|P|$, we have
\begin{equation}\label{|S|}
|S|=2|P|-|R|.
\end{equation}
We have
$$
|P|=\sum_{lm\leq X}  \ \sum_{\stackrel{r,s>0 \atop lr+ms=N}{ms\leq M }} 1
=\sum_{{lm\leq X}\atop{m\leq M}}  \ \sum_{\stackrel{ r,s>0 \atop lr+ms=N}{ s\leq \frac{M}{m}}} 1.
$$
Note that we can drop the condition $m\leq M$ in the outer sum, because
if $m>M$, then $s$ satisfies $0<s\leq M/m<1$, resulting an empty sum. Therefore,
\begin{equation}\label{P}
|P|=\sum_{lm\leq X}  \sum_{\stackrel{r,s>0 \atop  lr+ms=N}{s\leq \frac{M}{m}}} 1 =\sum_{lm\leq X}A(l,m),
\end{equation}
where $A(l,m)$ counts the positive integer pairs $(r,s)$, which are the solutions of $lr+ms=N$ such that
$0<s\leq M/m$. Let $d=\gcd(l,m)$ and write $l=d\lambda$, $m=d\mu$. If $d \nmid N$, then $A(l,m)=0$.
If $d\mid N$, the general form of the solutions is as follows:
$r=r_0+\mu h$ and $s=s_0-\lambda h$,
where $(r_0, s_0)$ is a particular solution and $h$ is an arbitrary integer. Hence $A(l,m)$ is equal to the number of integers $h$ satisfying  
\begin{equation}\label{cond-h}
	r_0+\mu h>0 \ \text{ and } \ s_0-\lambda h>0,
\end{equation}
together with
\begin{equation}\label{cond-s}
	0<s_0-\lambda h\leq \frac{M}{m}.
\end{equation}
From \eqref{cond-h}, we have
\begin{equation}\label{cond-h-1}
	-\frac{r_0}{\mu}< h< \frac{s_0}{\lambda},
\end{equation}
and using \eqref{cond-s}, we have
\begin{equation}\label{cond-h-2}
	\frac{s_0}{\lambda}-\frac{M}{m\lambda}\leq  h< \frac{s_0}{\lambda}.
\end{equation}
So, $A(l,m)$ denotes the number of integers $h$ satisfying \eqref{cond-h-1} and \eqref{cond-h-2} simultaneously. Hence
we count the number of $h$ such that
\begin{eqnarray*}
	\max\left(-\frac{r_0}{\mu},\frac{s_0}{\lambda}-\frac{M}{m\lambda}\right)\leq h< \frac{s_0}{\lambda}.
\end{eqnarray*}
Note that since $lr_0+ ms_0=N>M$, we have $\displaystyle{M<ms_0+\frac{r_0 d\lambda \mu}{\mu}}$, hence
$\displaystyle{-\frac{r_0}{\mu}<\frac{s_0}{\lambda}-\frac{M}{m\lambda}}$.
Therefore,
$$
A(l,m)=\left[\frac{s_0}{\lambda}\right]-\left[\frac{s_0}{\lambda}-\frac{M}{m\lambda}\right]
=\frac{M}{\lambda m}+ O(1) =\frac{M}{d\lambda \mu}+ O(1).
$$
Substituting this in \eqref{P}, we obtain
\begin{equation*}
|P|=M\sum_{d|N \atop d\leq \sqrt{X}}\frac{1}{d}\sum_{\lambda\mu \leq {X}/{d^2} \atop (\lambda,\mu)=1}\frac{1}{\lambda\mu}+O(X\log X)	
=M\sum_{d|N \atop d\leq \sqrt{X}}\frac{1}{d} \ \tau\left(\frac{X}{d^2}\right)+O(X\log X)	,
\end{equation*}
where $\displaystyle{\tau(y)=\sum_{\lambda\mu \leq y \atop (\lambda,\mu)=1}\frac{1}{\lambda\mu}}$.
From \cite[p. 203-204]{AEI}, we have (for sufficiently large $y$),
$$
\tau(y)=\frac{3}{\pi^2}\log^2y+O(\log y).
$$
So,
\begin{align*}
M\sum_{d|N \atop d\leq \sqrt{X}}\frac{1}{d} \ \tau\left(\frac{X}{d^2}\right)
&=\frac{3}{\pi^2}M\sum_{d|N\atop d\leq \sqrt{X} }\frac{1}{d}\log^2\left(\frac{X}{d^2}\right)+O\left(M\sum_{d|N\atop d\leq \sqrt{X} }\frac{1}{d}\log \left(\frac{X}{d^2}\right)\right)\\
&=\frac{3}{\pi^2}M\sum_{d|N \atop d\leq \sqrt{X}}\frac{1}{d}\left(\log^2 X+4\log^2 d-4(\log d)(\log X)\right)+O\left(M\sigma_{-1}(N)\log X\right),
\end{align*}
and thus
\begin{align*}
&M\sum_{d|N \atop d\leq \sqrt{X}}\frac{1}{d} \ \tau\left(\frac{X}{d^2}\right)\\
	&=\frac{3}{\pi^2}M(\log X)^2 \lambda_0(N) - \frac{12}{\pi^2}M (\log X) \lambda_1(N) +\frac{12}{\pi^2}M\lambda_2(N) +O\left(M\sigma_{-1}(N)\log X\right),
\end{align*}
where $\displaystyle{\lambda_0(N)=\sum_{d|N \atop d\leq \sqrt{X}}\frac{1}{d}, \lambda_1(N)=\sum_{d|N \atop d\leq \sqrt{X}}\frac{\log d}{d}}$ and $\displaystyle{\lambda_2(N)=\sum_{d|N \atop d\leq \sqrt{X}}\frac{\log^2 d}{d}}$.

Note that for $i=0,1,2$, in the sum $\lambda_i(N)$, we can remove the restriction $d \le \sqrt X$ as the resulting difference is bounded (in absolute value)
by $d(N) \log^2 N/\sqrt X$. Since $X \ge \sqrt{N-1}$, this is negligible and we have
\begin{align*}
&M\sum_{d|N \atop d\leq \sqrt{X}}\frac{1}{d} \ \tau\left(\frac{X}{d^2}\right)\\
&=\frac{3}{\pi^2}M(\log X)^2 \sum_{d|N \atop}\frac{1}{d} - \frac{12}{\pi^2}M \log X \sum_{d|N \atop}\frac{\log d}{d}
+\frac{12}{\pi^2}M \sum_{d|N \atop}\frac{\log^2 d}{d} +O\left(M\sigma_{-1}(N)\log X\right).
\end{align*}
This is clearly indicative of the terms to expect for an extension of Estermann's work \cite{TE}  in this aspect, which we do not take up here.
It follows from \cite[p. 208]{AEI} that
$$
\sum_{d|N \atop}\frac{\log d}{d}=O(\sigma_{-1}(N) \log \log N).
$$
Thus for $M$ sufficiently large such that $1\le M\le N/2$, we get
\begin{equation}\label{|P|-final}
M\sum_{d|N \atop d\leq \sqrt{X}}\frac{1}{d} \ \tau\left(\frac{X}{d^2}\right)=\frac{3}{\pi^2}M\sigma_{-1}(N) \log^2 X + O\left(M\sigma_{-1}(N)\log N \log\log N\right).
\end{equation}

For an asymptotic estimate of $|P|$, we need $O(X\log X)$ to be $o(M\sigma_{-1}(N)\log^2 X)$. This is true when 
\begin{equation}\label{o-condition}
\frac{N}{(\sigma_{-1}(N)\log N)^2}=o(M).
\end{equation}
To see this, we write $X=\sqrt{M(N-M)}= M\sqrt{\frac{N}{M}-1}$. From \eqref{o-condition}, we get $\left({N}/{M}\right)^{1/2} = o(\sigma_{-1}(N)\log N),$ and $X=o(M\sigma_{-1}(N)\log N).$ Hence $O(X\log X)$ is
$o(M\sigma_{-1}(N)\log^2 X)$ for all $M$ satisfying \eqref{o-condition}.

Moreover, for ${N}/{(\sigma_{-1}(N)\log\log N)^2}\le M\le {N}/{2}$, we have
\begin{equation}\label{ex-add-divisor-1}
|P| =\frac{3}{\pi^2}M\sigma_{-1}(N)\log ^2X +O(M\sigma_{-1}(N)\log N \log\log N),
\end{equation}
and for $M\le {N}/{(\sigma_{-1}(N)\log\log N)^2}$ such that ${N}/{(\sigma_{-1}(N)\log N)^2}=o(M)$, we have
\begin{equation}\label{ex-add-divisor-2}
|P| =\frac{3}{\pi^2}M\sigma_{-1}(N)\log ^2X +O(X\log X).
\end{equation}

Now it remains to estimate $|R|$, for which we write
$$
|R|\le \sum_{lm\leq X \atop ms \leq M}B(l,m),
$$
where $B(l,m)$ is the number of solutions ($r,s$) of $lr+ms=N$ such that $r,s>0$ and $rs\leq X$.
As before, let $d=\gcd(l,m)$ and write $l=d\lambda$, $m=d\mu$.
Note that if $d \nmid N$, then $B(l,m)=0$. If $d\mid N$, it is equal to the number of integers $h$ satisfying 
\begin{equation}\label{cond-h-R}
r_0+\mu h>0 \ \text{ and } \ s_0-\lambda h>0
\end{equation}
such that
\begin{equation}\label{cond-R}
(r_0+\mu h)(s_0-\lambda h)\leq X,
\end{equation}
where $(r_0, s_0)$ is a particular solution. 
Now since $l(r_0+\mu h)+m(s_0-\lambda h)=N$, we can write 
$$r_0+\mu h=\frac{N}{l}u \text{ and } s_0-\lambda h=\frac{N}{m}(1-u),
$$
and the above conditions are equivalent to the conditions
$$
0<u<1,
$$
and
\begin{equation}\label{quad}
-u^2+u-\frac{Xlm}{N^2}\leq 0.
\end{equation}
Now, we consider the quadratic equation 
$$
-u^2+u-\frac{Xlm}{N^2}= 0.
$$
Note that its discriminant $\Delta=1-4Xlm/N^2$. Since $lm \le X$ and $X \le N/2$, we get $\Delta \ge 0$.
Hence we can rewrite \eqref{quad} as
$$
\left(u-\left(\frac{1}{2}+j\right)\right)\left(\left(\frac{1}{2}-j\right)-u\right)\leq 0,
$$
where 
$$
 j:=\frac{1}{2}\sqrt{1-\frac{4Xlm}{N^2}}.
 $$
 Note that $j \le 1/2$. The above inequality says that exactly one of the inequalities 
$$
u-\left(\frac{1}{2}+j\right)\geq 0, \ \left(\frac{1}{2}-j\right)-u\geq 0
$$ is true.
Since $0<u<1$, we have either
\begin{equation}\label{u-1}
0<u\leq \frac{1}{2}-j,
\end{equation}
or
\begin{equation}\label{u-2}
\frac{1}{2}+j\leq u<1.
\end{equation}
From \eqref{u-1}, we have $0<r_0+\mu h\le \frac{N}{l}\left(\frac{1}{2}-j\right)$. Therefore,
\begin{eqnarray}
-\frac{r_0}{\mu}< h\leq \frac{N}{l\mu}\left(\frac{1}{2}-j\right)-\frac{r_0}{\mu}.
\end{eqnarray}
Note that the length of the half-open interval $\left(-\frac{r_0}{\mu}, \frac{N}{l\mu}(\frac{1}{2}-j)-\frac{r_0}{\mu}\right]$ is
$\frac{N}{l\mu}(\frac{1}{2}-j)$.
Similarly, from \eqref{u-2}, we have
$$
\frac{N}{l}\left(\frac{1}{2}+j\right)<r_0+\mu h\leq \frac{N}{l}.
$$
Therefore,
\begin{eqnarray}
	\frac{N}{l\mu}\left(\frac{1}{2}+j\right)-\frac{r_0}{\mu}\leq h< \frac{N}{l\mu}-\frac{r_0}{\mu}.
\end{eqnarray}
So, the length of the half-open interval $\left[ \frac{N}{l\mu}\left(\frac{1}{2}+j\right)-\frac{r_0}{\mu},\frac{N}{l\mu}-\frac{r_0}{\mu}\right)$ is also
$\frac{N}{l\mu}(\frac{1}{2}-j)$. Now
$$
\frac{N}{l\mu}\left(\frac{1}{2}-j\right) 
=\frac{Nd}{2lm}(1-2j) 
=\frac{Nd}{2lm}\frac{(1-4j^2)}{(1+2j)} \le \frac{Nd}{2lm}(1-4j^2)
= \frac{Nd}{2lm}\cdot\frac{4lmX}{N^2} 
=\frac{2dX}{N}.
$$
So, the admissible values of $h$ are thus confined to two (half-open) intervals each of which is of length at most $2dX/N$.
Hence $B(l,m) \le {4dX}/{N}$, where $d=\gcd(l,m)$ is a divisor of $N$. Thus
$$
	|R| \leq \sum_{lm\leq X}B(l,m) \nonumber
	\leq \sum_{d|N \atop d< \sqrt{X}}\sum_{\lambda\mu \leq{X}/{d^2} \atop (\lambda,\mu)=1}\frac{4dX}{N} \nonumber
	\leq \frac{4X}{N}\sum_{d|N }d\sum_{\lambda\mu \leq{X}/{d^2} }1 \nonumber
	= \frac{4X}{N}\sum_{d|N }d\sum_{\nu \leq{X}/{d^2} }d(\nu) \nonumber.
$$
Using the average order of the divisor function, we get
$$
|R| = O\left(\frac{X}{N}\sum_{d|N }\frac{X}{d}\log \left(\frac{X}{d^2}\right)\right) \nonumber
= O\left(\frac{X^2}{N}  \log X\sum_{d|N }\frac{1}{d}\right)=O\left(\frac{X^2}{N} \sigma_{-1}(N) \log X\right).
$$
Note that $X^2/N \le M(N-M)/N \le M$. Hence
\begin{equation}\label{|R|-final}
|R|=O\left(M \sigma_{-1}(N)\log X \right).
\end{equation}
This completes the proof of Theorem \ref{thm-add-divisor}.

\section{Proof of Theorem \ref{thm-add-divisor-2}}
As in \S1.1, we have
\begin{align}\label{M-small}
\sum_{n \le M}d(n)d(N-n)
= 2 \sum_{t \le \sqrt{N}} \sum_{n \le M \atop n \equiv N \bmod t} d(n) +O(M \log M),
\end{align}
using \eqref{small-part}. To prove \eqref{small-part}, we use the following estimate of Linnik and Vinogradov \cite{LV}:
for $0< \alpha<1/2$, we have as $x \to \infty$,
$$
\sum_{n \le x \atop n \equiv a \bmod k} d(n) \ll \frac{\phi(k)}{k^2} x \log x,
$$
uniformly for $a,k$ with $\gcd(a,k)=1$, provided $k < x^{1-\alpha}$.
Here $\phi$ denotes the Euler $\phi$-function.
Now denoting $\delta=\gcd(t,N)$, we write
\begin{align*}
\sum_{n \le M}d(n)\sum_{t|N-n \atop \sqrt{N-n} \le t \le \sqrt{N}} 1 & \le \sum_{\sqrt{N-M} < t \le \sqrt{N}}\sum_{n_1\le {M}/{\delta}\atop n_1\equiv {N}/{\delta} \bmod {t}/{\delta}}d(n_1\delta)\\
& \le \sum_{\sqrt{N-M} < t \le \sqrt{N}}d((N,t))\sum_{n_1\le {M}/{\delta}\atop n_1\equiv {N}/{\delta} \bmod {t}/{\delta}}d(n_1).
\end{align*}
Applying the above estimate of Linnik and Vinogradov, we get
\begin{align*}
\sum_{n \le M}d(n)\sum_{t|N-n \atop \sqrt{N-n} \le t \le \sqrt{N}} 1
&\ll M\log M\left(\sum_{\sqrt{N-M} < t \le \sqrt{N}}\frac{d(t)}{t}\right).
\end{align*}
In the sum on the right hand side, there are about $\sqrt{N}-\sqrt{N-M} = {M}/{(\sqrt{N}+\sqrt{N-M})}$ many summands and every summand is $\ll {N^{\epsilon}}/{\sqrt{N-M}}$ for any $\epsilon>0$. This proves \eqref{small-part}.


Now we study the sum of the right hand side of \eqref{M-small}. Using \eqref{eq-PV}, we get
\begin{align*}
2\sum_{t \le \sqrt{N}} \sum_{n \le M \atop n \equiv N \bmod t} d(n) =& \ 2\sum_{t \le \sqrt{N}}\left\{ \frac{M}{t} \sum_{r \mid t} \frac{c_r(N)}{r} \left( \log\frac{M}{r^2} + 2\gamma-1 \right) + O((M^{1/3}+t^{1/2})M^{\epsilon}) \right\}\\
=&\ 2M\log M\sum_{t\le\sqrt{N}}\frac{1}{t}\sum_{r \mid t} \frac{c_r(N)}{r}-4M\sum_{t\le\sqrt{N}}\frac{1}{t}\sum_{r \mid t} \frac{c_r(N)\log r}{r}\\
& + 2M(2\gamma-1)\sum_{t\le\sqrt{N}}\frac{1}{t}\sum_{r \mid t} \frac{c_r(N)}{r}+O(M).
\end{align*}
Now we have
\begin{align*}
\sum_{t\le\sqrt{N}}\frac{1}{t}\sum_{r \mid t} \frac{c_r(N)}{r}
=\sum_{r\le\sqrt{N}} \frac{c_r(N)}{r} \sum_{t\le\sqrt{N}\atop r|t}\frac{1}{t}
&=\sum_{r\le\sqrt{N}}\frac{c_r(N)}{r^2}\left(\log\left(\frac{\sqrt{N}}{r}\right)+O(1)\right)\\
&=\frac{\log N}{2}\mathfrak S_0(N) - \mathfrak S_1(N) + O_N(1),
\end{align*}
as from \eqref{sigma_s}, we know that for a fixed $N$, the series $\sum_{r\ge1}{c_r(N)}/{r^s}$ converges absolutely for $\Re(s)>1$.
Similarly,
\begin{align*}
\sum_{t\le\sqrt{N}}\frac{1}{t}\sum_{r \mid t} \frac{c_r(N)\log r}{r}
=\sum_{r\le\sqrt{N}} \frac{c_r(N)\log r}{r} \sum_{t\le\sqrt{N}\atop r|t}\frac{1}{t}
&=\sum_{r\le\sqrt{N}}\frac{c_r(N)\log r}{r^2}\left(\log\left(\frac{\sqrt{N}}{r}\right)+O(1)\right)\\
&=\frac{\log N}{2}\mathfrak S_1(N) - \mathfrak S_2(N) + O_N(\log M),
\end{align*}
where we can take the same implied constant as in the previous estimate. Using the last two estimates, we conclude the proof of \thmref{thm-add-divisor-2}.


\section{Proof of Theorem \ref{thm-add-conv}}
For the proof of Theorem \ref{thm-add-conv}, we need an orthogonality statement (see \cite[Lemma 9 (a)]{GMP}).
\begin{lem}\label{lem3}
For $M \le N$,
$$
\sum_{n< M}c_r(n)c_s(N-n)=\delta_{r,s}Mc_r(N)+O(rs\log (rs)),
$$
where $\delta_{r,s}$ denotes the Kronecker delta function.
\end{lem}

Following the proof of \cite[Theorem 2]{CMS}, let $U$ be a parameter tending to infinity which is to be chosen later and we write
$\sum_{n< M}f(n)g(N-n)=A+B,$
where
$$
A:=\sum_{n< M}\sum_{r,s \atop rs\leq U}a_{f}(r)a_{g}(s)c_r(n)c_s(N-n)
\ \text{ and } \
B:=\sum_{n< M}\sum_{r,s \atop rs>U}a_{f}(r)a_{g}(s)c_r(n)c_s(N-n).
$$
Now we use the growth conditions on $a_{f}(r)$ and $a_{g}(s)$ to interchange the summations and then Lemma \ref{lem3} to get
$$A=C+D+O\left(\sum_{r,s \atop rs\leq U}\frac{1}{(rs)^{1+\delta}}rs\log rs\right),$$
where $C=M\sum_{r \ge 1}a_{f}(r)a_{g}(r)c_r(N)$ and $D=-M\sum_{r^2>U}a_{f}(r)a_{g}(r)c_r(N)$.
It has been deduced in \cite[Section 4.1]{CMS} that
$$O\left(\sum_{r,s \atop rs\leq U}\frac{1}{(rs)^{1+\delta}}rs\log rs\right)
=\begin{cases}
O(U^{1-\delta}\log^2 U) & \text{if $\delta<1$},\\
O(\log^3U)& \text{if $\delta=1$}, \\
O(1)& \text{if $\delta>1$}.
\end{cases}$$
The growth conditions on $a_{f}(r)$ and $a_{g}(r)$ also imply that $D=O\left(M\sum_{r>\sqrt{U}}{c_r(N)}/{r^{2+2\delta}}\right)$.
Following \cite[Section 4.1]{CMS}, we also get
$$
\sum_{r>\sqrt{U}}{c_r(N)}/{r^{2+2\delta}}=O\left(\frac{e^{c\sqrt{\log N}}d(N)}{U^{1/2+\delta}}\right),
$$
for some constant $c>0$. Thus, we have 
$D \ll_N {M}/{U^{1/2+\delta}}$,
and therefore
$$
A=\begin{cases}
M\sum_{r\ge 1}a_{f}(r)a_{g}(r)c_r(N)+O\left(\frac{M}{U^{1/2+\delta}}\right)+O(U^{1-\delta}\log^2 U) & \text{if $\delta<1$},\\
M\sum_{r\ge 1}a_{f}(r)a_{g}(r)c_r(N)+O\left(\frac{M}{U^{3/2}}\right)+ O(\log^3U)& \text{if $\delta=1$}, \\
M\sum_{r\ge 1}a_{f}(r)a_{g}(r)c_r(N)+O\left(\frac{M}{U^{1/2+\delta}}\right)+ O(1)& \text{if $\delta>1$}.
\end{cases}
$$

The estimation for term $B$ is slightly different, as we need to use
\thmref{thm-add-divisor} and \corref{cor-add-divisor} here, in place of \eqref{shifted-divisor} as in \cite[Section 4.3]{CMS}.
Using hypotheses about $a_{f}(r)$ and $a_{g}(s)$, we get
\begin{align*}
	B &\ll \sum_{n< M}\sum_{r,s \atop rs>U}\frac{1}{(rs)^{1+\delta}}\left|\sum_{r'|r \atop r'|n}r'\mu(r/r')\sum_{s'|s \atop s'|N-n}s'\mu(s/s')\right|\\
	&= \sum_{r,s \atop rs>U}\frac{1}{(rs)^{1+\delta}}\left|\sum_{r'|r }r'\mu(r/r')\sum_{s'|s \atop }s'\mu(s/s')\sum_{\stackrel{n< M}{r'|n, \ s'|N-n}}1 \right|\\
	&= \sum_{r,s \atop rs>U}\frac{1}{(rs)^{1+\delta}}\left|\sum_{r'|r}r'\mu(r/r')\sum_{s'|s}s'\mu(s/s')\sum_{m< M/r' \atop r'm\equiv N \bmod {s'}}1\right|.
\end{align*}
Let $l=\gcd(r',s')$. Hence if we write $r'=lr''$ and $s'=ls''$, we get $\gcd(r'',s'')=1$. Also note that $l \mid N$. Thus
\begin{align*}
B&\ll \sum_{r,s \atop rs>U}\frac{1}{(rs)^{1+\delta}}\sum_{r'|r}r'\sum_{s'|s}s'\sum_{m< M/r' \atop r'm\equiv N \bmod {s'}}1 \\
&= \sum_{l \mid N}\sum_{\stackrel{r,s}{rs>U \atop l|r,l|s}}\frac{1}{(rs)^{1+\delta}}\sum_{r''|r/l}lr''\sum_{s''|s/l}ls''\sum_{m< M/lr'' \atop r''m\equiv (N/l) \bmod {s''}}1\ .
\end{align*}
Writing $r=lr_0$ and $s=ls_0$, we get
\begin{align*}
B&\ll \sum_{l \mid N}\frac{1}{l^{2\delta}}\sum_{r_0,s_0 \atop r_0s_0>U/l^2}\frac{1}{(r_0s_0)^{1+\delta}}\sum_{r''|r_0}r''\sum_{s''|s_0}s''
\sum_{m< M/lr'' \atop r''m\equiv (N/l) \bmod {s''}}1 = E+F,
\end{align*}
where 
$$E:=\sum_{l \mid N}\frac{1}{l^{2\delta}}\sum_{r_0,s_0 \atop r_0s_0>U/l^2}\frac{1}{(r_0s_0)^{1+\delta}}\sum_{r''|r_0}r''\sum_{s''|s_0 \atop s''\leq M/lr''}s''\sum_{m< M/lr'' \atop r''m\equiv (N/l) \bmod {s''}}1$$
and 
$$F:=\sum_{l \mid N}\frac{1}{l^{2\delta}}\sum_{r_0,s_0 \atop r_0s_0>U/l^2}\frac{1}{(r_0s_0)^{1+\delta}}\sum_{r''|r_0}r''\sum_{s''|s_0 \atop s''> M/lr''}s''\sum_{m< M/lr'' \atop r''m\equiv (N/l) \bmod {s''}}1\ .$$
Note that for $E$, we have $s''\leq M/lr''$. Since
$$
\sum_{m< M/lr'' \atop r''m\equiv (N/l) \bmod {s''}}1\leq \frac{M}{lr''s''}+1,
$$
following \cite[Section 4.3]{CMS}, we get
$$E\ll M\sum_{l \mid N}\frac{1}{l^{1+2\delta}}\sum_{r_0,s_0 \atop r_0s_0>U/l^2}\frac{d(r_0)d(s_0)}{(r_0s_0)^{1+\delta}}\ll \frac{M\log^4U}{U^\delta}.$$

Following the treatment of the term $F$ as in \cite[Section 4.3]{CMS}, here we get
\begin{align*}
	F \ll \frac{1}{M^\delta}\sum_{l \mid N}\frac{1}{l^\delta}\sum_{r'' }\sum_{s'' }\sum_{m< M/lr'' \atop r''m\equiv (N/l) \bmod {s''}}1\ .
\end{align*}
Now the sum
$$
	\sum_{r'' \ge 1}\sum_{s'' \ge 1}\sum_{m< M/lr'' \atop r''m\equiv (N/l) \bmod {s''}}1
	=\sum_{r'' \ge 1}\sum_{s'' \ge 1}\sum_{mr''< M/l \atop s''|r''m-N/l}1 
	=\sum_{r'' \ge 1}\sum_{s'' \ge 1}\sum_{\stackrel{n< M/l}{r''|n \atop s''|N/l-n}}1 
	=\sum_{n< M/l}d(n)d(N/l-n) .
$$
So, finally we have
\begin{equation}
	F\ll \frac{1}{M^\delta}\sum_{l \mid N}\frac{1}{l^\delta}\sum_{n< M/l}d(n)d(N/l-n).
\end{equation}

First we show that the additive convolution sub-sum $\sum_{n\le M}d(n)d(N-n)\ll M\log^3N$ when $M=o(N)$. Applying the Cauchy-Schwarz inequality we write , 
\begin{align*}
\sum_{n\le M}d(n)d(N-n) &\le \left(\sum_{n\le M}d(n)^2\right)^{{1}/{2}}\left(\sum_{n\le M}d(N-n)^2\right)^{{1}/{2}}
\end{align*}
We apply an asymptotic formula of the moments of the divisor function, due to Wilson \cite{BMW}.
\begin{lem}\label{BMW} 
For an integer $r\ge 1$, 
$$
\sum_{n\le x}d(n)^r \sim C_r x\left(\log x \right)^{2^r-1},
$$
as $x \to \infty$, where
$$
C_r = \frac{1}{(2^r-1)!}\prod_{p \ \text{prime}} \left\{\left(1-\frac{1}{p}\right)^{2^r}\left(\sum_{\alpha\ge0}\frac{(\alpha+1)^r}{p^{\alpha}}\right)\right\}.
$$   
\end{lem}
Applying \lemref{BMW}, we have
$$
\sum_{n\le M} d(n)^2 \ll M \log^3 M,
$$
and
\begin{align*}
\sum_{n\le M}d(N-n)^2&= \sum_{n< N}d(n)^2-\sum_{n< N-M}d(n)^2\sim C_2N\log^3N - C_2(N-M)\log^3(N-M),
\end{align*}
which is also $\ll M\log^3N$, as $M=o(N)$. Therefore, for all $M=o(N)$, we have
\begin{align*}\label{F}
	F\ll \frac{1}{M^\delta}\sum_{l \mid N}\frac{1}{l^\delta}\left(\frac{M}{l}\log^3\left(\frac{N}{l}\right)\right)
    \ll M^{1-\delta}\log^3N\sum_{l|N}\frac{1}{l^{1+\delta}}
    \ll_N M^{1-\delta}.
\end{align*}
To optimize all the error terms, we choose
$$U=\begin{cases}
	M\log^2M & \text{if $\delta<1$},\\
	M\log M& \text{if $\delta=1$},\\
	M^{(1/\delta)}(\log M)^{4/\delta}& \text{if $\delta>1$}.
\end{cases}$$
These choices yield
$$\sum_{n< M}f(n)g(N-n)=\begin{cases}
	M\sum_{r\ge 1}a_{f}(r)a_{g}(r)c_r(N)+O(M^{1-\delta}(\log M)^{4-2\delta}) & \text{if $\delta<1$},\\
	M\sum_{r\ge 1}a_{f}(r)a_{g}(r)c_r(N)+ O(\log^3M)& \text{if $\delta=1$}, \\
	M\sum_{r\ge 1}a_{f}(r)a_{g}(r)c_r(N)+ O(1)& \text{if $\delta>1$}.
\end{cases}$$
This concludes the proof of Theorem \ref{thm-add-conv}.

\begin{rmk}\label{rmk-2}\rm
The proof suggests that the implied constants in the error terms depend on $N$ and have order
bounded by $e^{c\sqrt{\log N}}d(N)$ for some constant $c >0$. If $N^\alpha \ll M$ for some $\alpha>0$ and $M=o(N)$, the
implied constants in the error terms can be made independent of $N$ by adjusting a $\log$ factor
as follows:
$$
\sum_{n< M}f(n)g(N-n)=\begin{cases}
M\sum_{r\ge 1}a_{f}(r)a_{g}(r)c_r(N)+O(M^{1-\delta}(\log M)^{4-\delta}) & \text{if $\delta<1$},\\
M\sum_{r\ge 1}a_{f}(r)a_{g}(r)c_r(N)+ O(\log^3M)& \text{if $\delta=1$}, \\
M\sum_{r\ge 1}a_{f}(r)a_{g}(r)c_r(N)+ O(1)& \text{if $\delta>1$}.
\end{cases}
$$
\end{rmk}

\subsection{Applications}
By the virtue of Theorem \ref{thm-add-conv}, we can now naturally extend and improve the corollaries that were obtained in \cite{GMP}.
Using the Ramanujan expansion \eqref{sigma_s}, we have the following corollary.
\begin{cor}\label{cor-sigma/n}
For $\alpha,\beta>0$, let $\delta:=\min\{ \alpha,\beta \}$. Then for $M<N$ sufficiently large and $M=o(N)$,
$$\mathlarger{\sum}_{n < M} \frac{\sigma_{\alpha}(n)}{n^{\alpha}}\frac{\sigma_{\beta}(N-n)}{(N-n)^{\beta}}=\begin{cases}
	\frac{M\zeta(\alpha+1)\zeta(\beta+1)}{\zeta(\alpha+\beta+2)}\frac{\sigma_{\alpha+\beta+1}(N)}{N^{\alpha+\beta+1}}+O(M^{1-\delta}(\log M)^{4-2\delta}) & \text{if $\delta<1$},\\
	\frac{M\zeta(\alpha+1)\zeta(\beta+1)}{\zeta(\alpha+\beta+2)}\frac{\sigma_{\alpha+\beta+1}(N)}{N^{\alpha+\beta+1}}+ O(\log^3M)& \text{if $\delta=1$}, \\
	\frac{M\zeta(\alpha+1)\zeta(\beta+1)}{\zeta(\alpha+\beta+2)}\frac{\sigma_{\alpha+\beta+1}(N)}{N^{\alpha+\beta+1}}+ O(1)& \text{if $\delta>1$}.
\end{cases}$$ 
where the implied constants may depend on $N$.
\end{cor}
This has been obtained in \cite[Theorem 12]{GMP} only in an asymptotic form and for $\alpha,\beta >1/2$.



\section*{Acknowledgement}
Research of B. Saha was supported by SERB grant SRG/2022/001011 and IIT Delhi IRD project MI02455.
Research of A. Sharma was sponsored by IIT Delhi IRD project MI02455. The authors are also thankful to
Prof. M. Ram Murty for his comments on an earlier version of this article.

\end{document}